\input btxmac.tex


\magnification=\magstep1 
\hsize=15truecm 
\vsize=22truecm
\hoffset=.65truecm 
\voffset=.3truecm


\hyphenation{Pro-po-si-tion}


\font\tenbb=msbm10     
\font\sevenbb=msbm7
\newfam\bbfam \def\bb{\fam\bbfam\tenbb}  
\textfont\bbfam=\tenbb
\scriptfont\bbfam=\sevenbb

\font\tenam=msam10     
\font\sevenam=msam7
\newfam\amfam   
\textfont\amfam=\tenam
\scriptfont\amfam=\sevenam

\font\bn=cmb10

\font\titlefont=cmr12 at 14pt
\let\authorfont=\rm
\let\titlerm=\elevenrm

\let\dinky=\sevenrm


\def\glq{\hbox{,\kern-.7pt,\kern.3pt}}  

\let\em=\it

\font\sevensl=cmsl8 at 7pt
\scriptfont\slfam=\sevensl

\let\titlefont=\twelverm  
\let\authorfont=\tenrm  
\let\sectheadfont=\titleit  


\def\secthead#1#2{\removelastskip\bigskip\goodbreak%
{\titlerm #1\enspace}{\sectheadfont #2}%
\nobreak\medskip\noindent}


\thickmuskip=4.5mu plus 4mu minus 3mu  
\mathsurround=0.5pt

\pretolerance=200  \hyphenpenalty=1000 \exhyphenpenalty=5000

\baselineskip=14pt \parskip=1.5pt \parindent=1.5em
\smallskipamount=4pt plus2pt minus1pt 
\medskipamount=7pt plus3pt minus2pt 
\bigskipamount=15pt plus5pt minus4pt 


\let\originalsqrt\sqrt
\def\varsqrt #1#2{{\setbox0=\hbox{$#1\originalsqrt{#2\hskip2pt}$}%
     \dimen0=\ht0 \advance\dimen0-.5pt
     \dimen2=\dimen0
     \advance\dimen2 -1.5pt
     \ifdim \dimen2 >10pt \advance\dimen2-.5pt \fi
     \ifdim \dimen2 >15pt \advance\dimen2 -1pt \fi
     \box0 \kern-.4pt
     \vrule width.4pt height\dimen0 depth -\dimen2 \kern.8pt
  }}


\chardef\at="40  

\def\c{{\rm c}}
\def\d{{\rm d}}
\def\Int{\hbox{\rm Int}}
\def\Lset{{\cal L}}
\let\N=\natn

\def\intz{{\bb Z}}     
\let\Z=\intz
\let\Q=\ratq

\def\[{{\rm[}} \def\]{\/{\rm]}}  
\def\({{\rm(}} \def\){\/{\rm)}}  
\def\divides{\mathrel\big\vert} 

\let\congr=\equiv    
\mathchardef\ideal="1945  
\mathchardef\propersubset="1824  
\def\subsetneq{\mathrel\propersubset}
\mathchardef\eop="1903  

\newcount\codenumber \codenumber=`:
\loop\relax\if\codenumber>"FFF \advance\codenumber\by -"1000\repeat
\advance\codenumber by "3000
\mathcode`:=\codenumber


\long\def\profess#1#2\endprofess {\ifdim\lastskip<\medskipamount%
  \removelastskip\medskip\penalty-500\fi
  \noindent{\bn#1.\ }{\sl#2\par}%
  \ifdim\lastskip<\smallskipamount\removelastskip\penalty-55\smallskip\fi}

\long\def\rmprofess#1#2\endrmprofess {\ifdim\lastskip<\medskipamount%
  \removelastskip\medskip\penalty-350\fi
  \noindent{\bn#1.\ }{\rm#2\par}%
  \ifdim\lastskip<\smallskipamount\removelastskip\penalty-100\smallskip\fi}

\def\proof{\noindent {\it Proof.\enspace}}
\def\endproof{\ $\eop$\par%
   \ifdim\lastskip<\medskipamount
   \removelastskip\penalty-50\medskip\fi}


\def\firstpagehead{\tenrm\  To appear in Monatsh.~Math.\hfill }
\def\runninghead{\hfill 
{\dinky  \ } \hfill }
\headline={\ifnum\pageno>1\runninghead \else\firstpagehead \fi}

\output={
          \shipout\vbox{\makeheadline
                        \vskip.1truein
                        \pagebody
                        \makefootline }
          \advancepageno
          \ifnum\outputpenalty>-20000 \else\dosupereject\fi
 }

\long\def\comment#1\endcomment{\relax}

\def\uniqueexprem{1}
\def\atomicrem{2}
\def\wellknownrem{3}
\def\factorizationlem{4}
\def\residuelem{5}
\def\irredsamefdlem{6}
\def\nplustwoex{7}
\def\nandmex{8}
\def\setoflengththm{9}
\def\noblockmonoidthm{10}

\centerline{\titlefont A construction of integer-valued polynomials}
\centerline{\titlefont  with prescribed sets of lengths of factorizations}
\bigskip
\centerline{\authorfont Sophie Frisch}
\bigskip
\noindent {\bf Abstract.}
For an arbitrary finite set $S$ of natural numbers greater $1$, we
construct $f\in\Int(\Z)=\{g\in\Q[x]\mid g(\Z)\subseteq \Z\}$ whose
set of lengths is $S$. The set of lengths of $f$ is the set of all
$n$ such that $f$ has a factorization as a product of $n$ irreducibles
in $\Int(\Z)$.
\bigskip
\secthead{1.}{Introduction}%
Non-unique factorization has long been studied in rings of
integers of number fields, see the monograph of Geroldinger and
Halter-Koch \cite{GerHKNUF06}. More recently, non-unique factorization 
in rings of polynomials has attracted attention, for instance
in $\Z_{p^n}[x]$, cf.~\cite{FreFriNUF11}, and in the ring of
integer-valued polynomials
$\Int(\Z)=\{g\in\Q[x]\mid g(\Z)\subseteq \Z\}$
(and its generalizations) \cite{CahChaEIVP95,ChapMcCFE05}.

We show that every finite set of natural numbers greater~$1$
occurs as the set of lengths of factorizations of an element of
$\Int(\Z)$ (Theorem \setoflengththm\ in section 4).

Our proof is constructive, and allows multiplicities of lengths
of factorizations to be specified. For example, given the
multiset $\{2,2,2,5,5\}$, we construct a polynomial that has three
different factorizations into $2$ irreducibles and two different
factorizations into $5$ irreducibles, and no other factorizations.
Perhaps a quick review of the vocabulary of factorizations is in order:

\profess{Notation and Conventions}
$R$ denotes a commutative ring with identity.
An element $r\in R$ is called {\em irreducible} in $R$ if
$r$ is a non-zero non-unit such that $r=ab$ with $a,b\in R$ implies 
that $a$ or $b$ is a unit. A {\em factorization} of $r$ in $R$
is an expression $r=s_1\ldots s_n$ of $r$ as a product of irreducible
elements in $R$. The number $n$ of irreducible factors is called the
{\em length} of the factorization. The {\em set of lengths} $\Lset(r)$ of
$r\in R$ is the set of all natural numbers $n$ such that $r$ has a 
factorization of length $n$ in $R$.

$R$ is called {\em atomic} if every non-zero non-unit of $R$ has a 
factorization in $R$. 
%
If $R$ is atomic, then for every non-zero non-unit $r\in R$ the
{\em elasticity of $r$} is defined as
$$\rho(r)=\sup\{{m\over n}\mid m,n\in\Lset(r)\}$$
and the elasticity of $R$ is $\rho(R)=\sup_{r\in R'}(\rho(r))$,
where $R'$ is the set of non-zero non-units of $R$. An atomic
domain $R$ is called {\em fully elastic} if every rational number
greater than $1$ occurs as $\rho(r)$ for some non-zero
non-unit $r\in R$.

Two elements $r,s\in R$ are called {\em associated} in $R$ if there
exists a unit $u\in R$ such that $r=us$.
Two factorizations of the same element
$r=r_1\cdot\ldots\cdot r_m=s_1\cdot\ldots\cdot s_n$ are called 
{\em essentially the same} if $m=n$ and, after re-indexing the $s_i$,
$r_j$ is associated to $s_j$ for $1\le j\le m$. Otherwise, the
factorizations are called {\em essentially different}.
\endprofess
\secthead{2.}{Review of factorization of integer-valued polynomials}%
In this section we recall some elementary properties of $\Int(\Z)$ and
the fixed divisor $\d(f)$, to be found in \cite{CahChaEIVP95}, 
\cite{CaCh97ivp} and \cite{ChapMcCFE05}.
The reader familiar with integer-valued polynomials is encouraged to skip 
to section 3.

\profess{Definition}
For 
$f\in\Z[x]$, 
\item{\rm (i)} the content $\c(f)$ is the ideal
of $\Z$ generated by the coefficients of $f$,
\item{\rm (ii)} the fixed divisor $\d(f)$ is the ideal of $\Z$ 
generated by the image $f(\Z)$.

By abuse of notation we will identify the principal ideals
$\c(f)$ and $\d(f)$ with their non-negative generators. Thus,
for $f=\sum_{k=0}^n a_kx^k\in\Z[x]$, 
$$\c(f)=\gcd{(a_k\mid k=0,\ldots, n)}
\quad\hbox{\rm\ and\ }\quad
\d(f)=\gcd{(f(c)\mid c\in\intz)}.$$
A polynomial $f\in\Z[x]$ is called primitive if $\c(f)=1$.
\endprofess

Recall that a primitive polynomial $f\in\Z[x]$ is irreducible in
$\Z[x]$ if and only if it is irreducible in $\Q[x]$.
Similarly, $f\in\Z[x]$ with $\d(f)=1$ is irreducible in $\Z[x]$ 
if and only if it is irreducible in $\Int(\Z)$. 

We denote $p$-adic valuation by $v_p$. Almost everything that we
need to know about the fixed divisor follows immediately from the
fact that $$v_p(\d(f))=\min_{c\in\Z}(v_p(f(c))).$$ In particular,
it is easy to deduce that for any $f,g\in\Z[x]$,
$$\d(f)\d(g) \divides \d(fg).$$
Unlike $\c(f)$, which satisfies $\c(f)\c(g)=\c(fg)$,
$\d(f)$ is not multiplicative: $\d(f)\d(g)$ is in general a 
proper divisor of $\d(fg)$.

\profess{Remark \uniqueexprem}
\item{\rm (i)}
Every non-zero polynomial $f\in\Q[x]$ can be written in a unique way as
$$f(x)={{a g(x)}\over b}
\quad\hbox{\rm with}\quad
g\in\Z[x],\; c(g)=1,\quad a,b\in\N,\; \gcd(a,b)=1.$$
\item{\rm (ii)}
When expressed as in {\rm (i)}, $f$ is in $\Int(\Z)$ if and only if
$b$ divides $\d(g)$. 
\item{\rm (iii)} For  non-constant $f\in\Int(\Z)$ expressed as in 
{\rm (i)} to be irreducible in $\Int(\Z)$ it is necessary that $a=1$
and $b=\d(g)$.
\endprofess

\proof (i) and (ii) are easy.
Ad {\rm (iii)}. Note that the only units in $\Int(\Z)$ are
$\pm 1$. By (ii), $b$ divides $\d(g)$. Let $\d(g)=bc$.
Then $f$ factors
as $a\cdot c\cdot (g/bc)$, where $(g/bc)$ is non-constant and
$ac$ is a unit only if $a=c=1$.
\endproof

\profess{Remark \atomicrem}
\item{\rm (i)}
Every non-zero polynomial $f\in\Q[x]$ can be written in a unique way
\(up to the sign of $a$ and the signs and indexing of the $g_i$\) as
$$f(x)={a\over b}\prod_{i\in I} g_i(x),$$
with $g_i$ primitive and irreducible in $\Z[x]$ for $i\in I$ \(a finite
set\) and $a\in\Z$, $b\in\N$ with $\gcd(a,b)=1$.

\item{\rm (ii)} A non-constant polynomial $f\in\Int(\Z)$ expressed as
in {\rm (i)} is irreducible in $\Int(\Z)$ if and only if $a=\pm 1$, 
$b=\d(\prod_{i\in I} g_i)$, and there do not exist
$\emptyset\ne J\subsetneq I$ and $b_1,b_2\in\N$
with $b_1b_2=b$ and $b_1=\d(\prod_{i\in J} g_i)$,
$b_2=\d(\prod_{i\in I\setminus J} g_i)$.

\item{\rm (iii)} $\Int(\Z)$ is atomic.

\item{\rm (iv)} Every non-zero non-unit $f\in \Int(\Z)$ has only
finitely many factorizations into irreducibles in $\Int(\Z)$.
\endprofess

\proof
Ad {\rm (ii)}. If $f$ is irreducible, the conditions on $f$ follow from
Remark \uniqueexprem\ (ii) and (iii). Conversely, if the conditions
hold, what chance does $f$ have to be reducible?
By Remark \uniqueexprem\ (ii), we cannot factor out
a non-unit constant, because no proper multiple of $b$ divides
$\d(\prod_{i\in I} g_i)$. Any non-constant irreducible factor would,
by Remark \uniqueexprem\ (iii), be of the kind $(\prod_{i\in J} g_i)/b_1$
with $b_1=\d(\prod_{i\in J} g_i)$, and its co-factor would be
$(\prod_{i\in I\setminus J} g_i)/b_2$ with $b_1b_2=b$ and $b_2$
a divisor of $\d(\prod_{i\in I\setminus J} g_i)$. Also, $b_2$
could not be a proper divisor of $\d(\prod_{i\in I\setminus J} g_i)$,
because otherwise $b_1b_2=b$ would be a proper divisor of
$\prod_{i\in I} g_i$. So, the existence of a non-constant irreducible
factor would imply the existence of $J$ and $b_1,b_2$ of the kind we
have excluded.
 
Ad {\rm (iii)}. With $f(x)={{a g(x)}/ b}$, $g=\prod_{i\in I}g_i$ 
as in {\rm (i)}, $\d(g)=cb$ for some $c\in \N$, and 
$f(x)={{ac g(x)}/\d(g)}$ with $g(x)/\d(g)\in\Int(\Z)$.
We can factor $ac$ into irreducibles in $\Z$, which are also irreducible
in $\Int(\Z)$.
Either $g(x)/\d(g)$ is irreducible, or {\rm (ii)} gives an
expression as a product of two non-constant factors of smaller degree.
By iteration we arrive at a factorization of $g(x)/\d(g)$ into
irreducibles.

Ad {\rm (iv)}. Let $f\in\Int(\Z)=(ag(x)/b)$ with $g=\prod_{i\in I}g_i$
as in {\rm (i)}.
Then all factorizations of
$f$ are of the form, for some $c\in\N$ such that $bc$ divides $\d(g)$,
$$f= a_1\ldots a_n c_1\ldots c_m
\prod_{j=1}^k {{\prod_{i\in I_j}g_i}\over {d_j}},$$
where $a=a_1\ldots a_n$ and $c=c_1\ldots c_m$ are factorizations into
primes in $\Z$, $I=I_1\cup\ldots\cup I_k$ is a partition of $I$ into
non-empty sets, $d_1\ldots d_k= bc$, $d_j=\d(\prod_{i\in I_j}g_i)$.
There are only finitely many such expressions.
\endproof

\profess{Remark \wellknownrem}
\item{\rm (i)}
The {\em binomial polynomials}
$${x\choose n}={{x(x-1)\ldots(x-n+1)}\over{n!}}
\quad\hbox{\rm for}\quad n\ge 0$$
are a basis of $\Int(\Z)$ as a free $\Z$-module.
\item{\rm (ii)}
$n!f\in\Z[x]$ for every $f\in\Int(\Z)$ of degree at most $n$.
\item{\rm (iii)}
Let $f\in\Z[x]$ primitive, $\deg f=n$ and $p$ prime.  Then 
$$v_p(\d(f))\le \sum_{k\ge 1}\left[{{n}\over{p^k}}\right]=v_p(n!).$$
In particular, if $p$ divides $\d(f)$ then $p\le\deg f$.
\endprofess

\proof
Ad (i). The binomial polynomials are in $\Int(\Z)$ and they form a
$\Q$-basis of $\Q[x]$. If a polynomial in $\Int(\Z)$ is written as a 
$\Q$-linear combination of binomial polynomials then an easy induction
shows that the coefficients must be integers. (ii) follows from (i).

Ad (iii). Let $g=f/d(f)$. Then $g\in\Int(\Z)$ and $\d(f)\Z=(\Z[x]:_{\Z} g)$.
Since $n!\in (\Z[x]:_{\Z} g)$ by (ii),
$\d(f)$ divides $n!$.
\endproof
\vfill\eject

\secthead{3.}{Useful Lemmata}%
\profess{Lemma \factorizationlem}
Let $p$ be a prime, $I\ne \emptyset$ a finite set and for $i\in I$,
$f_i\in\Z[x]$ primitive and irreducible in $\Z[x]$ such that
$\d(\prod_{i\in I} f_i)=p$.  Let 
$$g(x)={{\prod_{i\in I} f_i}\over p}.$$
Then every factorization of $g$ in $\Int(\Z)$ is essentially the same
as one of the following:
$$g(x)={{\prod_{j\in J} f_j}\over p}\cdot \prod_{i\in I\setminus J}f_i, $$ 
where $J\subseteq I$ is minimal such that $\d(\prod_{i\in J} f_j)=p$.
\endprofess

\proof
Follows from Remark \uniqueexprem\ (iii) and the fact that 
$\d(f)\d(h)$ divides $\d(fh)$ for all $f,h\in\Z[x]$.
\endproof

The following two easy lemmata are constructive, since the Euclidean
algorithm makes the Chinese Remainder Theorem in $\Z$ effective.

\profess{Lemma \residuelem}
For every prime $p\in\Z$, we can construct a complete
system of residues mod $p$ that does not contain a complete
system of residues modulo any other prime.
\endprofess

\proof
By the Chinese Remainder Theorem we solve, for each $k=1,\ldots, p$
the system of congruences $s_k= k$ mod $p$ and $s_k=1$ mod $q$ for 
every prime $q<p$.
\endproof

\profess{Lemma \irredsamefdlem}
Given finitely many non-constant monic polynomials $f_i\in\Z[x]$,
$i\in I$, we can construct monic irreducible polynomials $F_i\in\Z[x]$,
pairwise non-associated in $\Q[x]$, with $\deg{F_i}=\deg{f_i}$, 
and with the following property:

Whenever we replace some of the $f_i$ by the corresponding $F_i$, setting
$g_i=F_i$ for $i\in J$ ($J$ an arbitrary subset of $I$) and $g_i=f_i$ for
$i\in I\setminus J$, then for all $K\subseteq I$,
$$\d(\prod_{i\in K}\kern-2pt g_i)= \d(\prod_{i\in K}\kern-2pt f_i).$$
\endprofess 

\proof
Let $n=\sum_{i\in I}\deg f_i$. Let $p_1,\ldots, p_s$ be all the primes
with $p_i\le n$, and set $\alpha_i=v_{p_i}(n!)$. Let $q>n$ be a prime.
For each $i\in I$,
we find by the Chinese Remainder Theorem the coefficients of 
a polynomial $\varphi_i\in (\prod_{k=1}^s p_k^{\alpha_k})\Z[x]$ of
smaller degree than $f_i$,
such that $F_i=f_i+\varphi_i$ satisfies Eisenstein's irreducibility
criterion with respect to the prime $q$. Then, with respect to some
linear ordering of $I$, if $F_i$ happens to be associated in $\Q[x]$
to any $F_j$ of smaller index, we add a suitable non-zero integer 
divisible by $q^2\prod_{k=1}^s p_k^{\alpha_k}$ to $F_i$, to make $F_i$
non-associated in $\Q[x]$ to all $F_j$ of smaller index.

The statement about the fixed divisor follows, because for every
$c\in\Z$ and every prime $p_i$ that could conceivably divide the 
fixed divisor, 
$$\prod_{i\in K}(g_i(c))\congr \prod_{i\in K}(f_i(c))
\quad\hbox{\rm mod}\; p_i^{\alpha_i},$$ 
where $p_i^{\alpha_i}$ is the highest power of $p_i$ that can
divide the fixed divisor of any monic polynomial of degree at
most $n$.
\endproof
\secthead{4.}{Constructing polynomials with prescribed sets of lengths}%
We precede the general construction by two illustrative examples of
special cases, corresponding to previous results by Cahen, Chabert,
Chapman and McClain.

\profess{Example \nplustwoex}
For every $n\ge0$, we can construct $H\in\Int(\Z)$ such that $H$ has
exactly two essentially different factorizations in $\Int(\Z)$,
one of length $2$ and one of length $n+2$.
\endprofess

\proof
Let $p>n+1$, $p$ prime. By Lemma \residuelem\ we construct a complete set
$a_1,\ldots,a_p$ of residues mod $p$ in $\Z$ that does not contain a
complete set of residues mod any prime $q<p$.
Let 
$$f(x)=(x-a_2)(x-a_3)\ldots (x-a_{p})\quad\hbox{\rm and}\quad
g(x)=(x-a_{n+2})(x-a_{n+3})\ldots (x-a_{p}).$$
By Lemma \irredsamefdlem, we construct monic irreducible polynomials
$F,G\in\Z[x]$, not associated in $\Q[x]$, with $\deg F=\deg f$,
$\deg G=\deg g$, such that any product of a selection of polynomials
from $(x-a_{1}),\ldots,(x-a_{n+1}), f(x), g(x)$ has the same fixed
divisor as the corresponding product with $f$ replaced by $F$ and
$g$ by $G$.

Let $$H(x)={{F(x)(x-a_{1})\ldots(x-a_{n+1})G(x)}\over p}.$$
By Lemma \factorizationlem,
$H$ factors into two irreducible polynomials in $\Int(\Z)$
$$H(x)=F(x)\cdot{{(x-a_{1})\ldots(x-a_{n+1})G(x)}\over p}$$
or into $n+2$ irreducible polynomials in $\Int(\Z)$ 
$$H(x)=
{{F(x)(x-a_{1})}\over p}\cdot (x-a_2)(x-a_3)\ldots (x-a_{n+1})G(x).$$
\endproof

\profess{Corollary} {\rm (Cahen and Chabert \cite{CahChaEIVP95}).}
$\rho\,(\Int(\Z))=\infty$.
\endprofess

\profess{Example \nandmex} For $1\le m\le n$, we can construct a polynomial
$H\in\Int(\Z)$ that has in $\Int(\Z)$ a factorization into $m+1$
irreducibles and an essentially different factorization into $n+1$
irreducibles, and no other essentially different factorization.
\endprofess

\proof
Let $p>mn$ be prime, $s=p-mn$. By Lemma \residuelem\ we construct a
complete system of residues $R$ mod $p$ that does not contain
a complete system of residues for any prime $q<p$. We index $R$ as
follows:
$$R=\{r(i,j)\mid 1\le i\le m,\> 1\le j\le n\}\cup \{b_1,\ldots, b_s\}.$$
Let $b(x)=\prod_{k=1}^s (x-b_k)$.
For $1\le i\le m$ let $f_i(x)=\prod_{k=1}^n (x-r(i,k))$ and for
$1\le j\le n$ let $g_j(x)=\prod_{k=1}^m (x-r(k,j))$.

By Lemma \irredsamefdlem, we construct monic irreducible polynomials 
$F_i, G_j\in \Z[x]$, pairwise non-associated in $\Q[x]$, such that
the product of any selection of the polynomials 
$(x-b_1),\ldots,(x-b_s), f_1,\ldots, f_m, g_1,\ldots, g_n$
has the same fixed divisor as the corresponding product in which
$f_i$ has been replaced by $F_i$ and $g_j$ by $G_j$ for $1\le i\le m$
and $1\le j\le n$.
Let $$H(x)={1\over p}b(x)\prod_{i=1}^m F_i(x)\prod_{j=1}^n G_j(x),$$
then, by Lemma \factorizationlem, $H$ has a factorization into $m+1$ 
irreducibles
$$H(x)=F_1(x)\cdot\ldots\cdot F_m(x)\cdot
{{b(x)G_1(x)\cdot\ldots\cdot G_n(x)}\over p}$$
and an essentially different factorization into $n+1$ irreducibles
$$H(x)={{b(x)F_1(x)\cdot\ldots\cdot F_m(x)}\over p}\cdot
G_1(x)\cdot\ldots\cdot G_n(x)$$
and no other essentially different factorization.
\endproof

\profess{Corollary} {\rm (Chapman and McClain \cite{ChapMcCFE05}).}
$\Int(\Z)$ is fully elastic.
\endprofess

\profess{Theorem \setoflengththm}
Given natural numbers $1\le m_1\le\ldots\le m_n$, we can construct
a polynomial $H\in\Int(\Z)$ that has exactly $n$ essentially different 
factorizations into irreducibles in $\Int(\Z)$, the lengths of these
factorizations being $m_1+1,\ldots, m_n+1$.
\endprofess

\proof 
Let $N=(\sum_{i=1}^n m_i)^2-\sum_{i=1}^n m_i^2$, and $p>N$ prime,
$s=p-N$.
By Lemma \residuelem, we construct a complete system of residues $R$
mod $p$ that does not contain a complete system of residues for any
prime $q<p$.
We partition $R$ into disjoint sets 
$R=R_0\cup \{t_1,\ldots, t_s\}$ with $\left|R_0\right|= N$.
The elements of $R_0$ are indexed as follows:
$$ R_0=\{r(k, h, i, j)\mid 1\le k\le n,\>
1\le h\le m_k,\> 1\le i\le n,\> 1\le j\le m_i;\> i\ne k\},$$
meaning we arrange the elements of $R_0$ in an $m\times m$ matrix
with $m=m_1+\ldots+m_n$, whose rows and columns are partitioned
into $n$ blocks of sizes $m_1,\ldots, m_n$. Now $r(k, h, i, j)$
designates the entry in the $h$-th row of the $k$-th block of rows
and the $j$-th column of the $i$-th block of columns.
Positions in the matrix whose row and column are each in block $i$
are left empty: there are no elements $r(k, h, i, j)$ with $i=k$.

For $1\le k\le n,$ $1\le h\le m_k$, let $S_{k,h}$ be the set of
entries in the $(k,h)$-th row:
$$S_{k,h}=\{r(k,h,i,j)\mid 1\le i\le n,\> i\ne k,\> 1\le j\le m_i\}.$$
For $1\le i\le n$, $1\le j\le m_i$, let $T_{i,j}$ be the set of
elements in the $(i,j)$-th column:
$$T_{i,j}=\{r(k,h,i,j)\mid 1\le k\le n,\> k\ne i,\> 1\le h\le m_k\}.$$

For $1\le k\le n$, $1\le h\le m_k$, set
$$f^{(k)}_h(x)=\prod_{r\in S_{k,h}}(x-r)\cdot\prod_{r\in T_{k,h}}(x-r).$$

Also, let $b(x)=\prod_{i=1}^s(x-t_i)$.

By Lemma \irredsamefdlem, we construct monic irreducible polynomials
$F^{(k)}_h$, pairwise non-associated in $\Q[x]$, with
$\deg F^{(k)}_h =\deg f^{(k)}_h$, such that any product of a selection
of polynomials from $(x-t_1),\ldots,(x-t_s)$ and $f^{(k)}_h$ for
$1\le k\le n$, $1\le h\le m_k$ has the same fixed divisor as the
corresponding product in which the $f^{(k)}_h$ have been replaced by
the $F^{(k)}_h$.
Let $$H(x)={1\over p}b(x)\prod_{k=1}^n\prod_{h=1}^{m_k} F^{(k)}_h(x).$$
Then $\deg H = N+p$; and
for each $i=1,\ldots, n$, $H$ has a factorization into $m_i+1$
irreducible polynomials in $\Int(\Z)$:
$$H(x)=F^{(i)}_1(x)\cdot\ldots\cdot F^{(i)}_{m_i}(x)\cdot
{{b(x)\prod_{k\ne i}\prod_{h=1}^{m_k}F^{(k)}_h(x)}\over p}$$
These factorizations are essentially different, since the $F^{(i)}_j$
are pairwise non-associated in $\Q[x]$ and hence in $\Int(\Z)$.

By Lemma \factorizationlem, $H$ has no further essentially different
factorizations.
This is so because a minimal subset with fixed divisor $p$ of the
polynomials $(x-t_i)$ for $1\le i\le s$ and $F^{(k)}_h$ for
$1\le k\le n$, $1\le h\le m_k$ must consist of all the linear
factors $(x-t_i)$ together with a minimal selection of $F^{(k)}_h$
such that all $r\in R_0$ occur as roots in the product of the
corresponding $f^{(k)}_h$.
For all linear factors $(x-r)$ with $r\in R_0$ to occur in a set of
polynomials $f^{(k)}_{h}$, it must contain for all but one $k$ all
$f^{(k)}_{h}$, $h=1,\ldots m_k$. If, for $i\ne k$, 
$f^{(k)}_{h}$ and $f^{(i)}_j$ are missing, then $r(k,h,i,j)$ and 
$r(i,j,k,h)$ do not occur among the roots of the polynomials $f^{(k)}_h$.
A set consisting of all $f^{(k)}_{h}$ for $n-1$ different values of $k$,
however, has the property that all linear factors $(x-r)$ for $r\in R_0$ 
occur.
\endproof

\profess{Corollary}
Every finite subset of $\N\setminus \{1\}$ occurs as the set of lengths
of a polynomial $f\in\Int(\Z)$.
\endprofess

\secthead{5.}{No transfer homomorphism to a block-monoid}%

For some monoids, results like the above Corollary have been shown
by means of transfer-homomorphisms to block monoids.
For instance, by Kainrath \cite{KaiFKIKG99}, in the case of a Krull
monoid with infinite class group such that every divisor class 
contains a prime divisor.

$\Int(\Z)$, however, doesn't admit this method:
We will show a property of the multiplicative monoid of 
$\Int(\Z)\setminus\{0\}$ that excludes the existence of
a transfer-homomorphism to a block monoid.

\profess{Theorem \noblockmonoidthm}
For every $n\ge 1$ there exist irreducible elements 
$H,G_1,\ldots,G_{n+1}$ in $\Int(\intz)$ such that
$xH(x)=G_1(x)\ldots G_{n+1}(x)$.
\endprofess

\proof
Let $p_1<p_2<\ldots<p_n$ be $n$ distinct odd primes, 
$P=\{p_1,p_2,\ldots,p_n\}$, and $Q$ the set of all primes $q\le p_n +n$.
By the Chinese remainder theorem construct 
$a_1,\ldots, a_n$ with $a_i\congr 0$ mod~$p_i$ and $a_i\congr 1$
mod~$q$ for all $q\in Q$ with $q\ne p_i$. Similarly, construct
$b_1,\ldots b_{p_n}$ such that, firstly, for all $p\in P$,
$b_k\congr k$ mod~$p$ if $k\le p$ and $b_k\congr 1$ mod~$p$ if $k>p$ 
and, secondly, $b_k\congr 1$ mod~$q$ for all $q\in Q\setminus P$.
So, for each $p_i\in P$, a complete set of residues mod $p_i$ is
given by $b_1,\ldots b_{p_i},a_i$, while all remaining $a_j$ and $b_k$
are congruent to $1$ mod~$p_i$. Also, all $a_j$ and $b_k$ are congruent
to $1$ for all primes in $Q\setminus P$.

Set $f(x)=(x-b_1)\ldots (x-b_{p_n})$ and let $F(x)$ be a monic
irreducible polynomial in $\intz[x]$ with $\deg F=\deg f$ such that
the fixed divisor of any product of a selection of polynomials from
$f(x),(x-a_1), \ldots, (x-a_n)$ is the same as the fixed divisor
of the corresponding set of polynomials in which $f$ has been replaced
by $F$. Such an $F$ exists by Lemma~\irredsamefdlem. Let
$$H(x)={{F(x)(x-a_1)\ldots(x-a_n)}\over{p_1\ldots p_n}}.$$
Then $H(x)$ is irreducible in $\Int(\intz)$, and
$$xH(x)={{xF(x)}\over{p_1\ldots p_n}}\cdot 
(x-a_1)\cdot\ldots\cdot (x-a_n),$$
where $xF(x)/(p_1\ldots p_n)$ and, of course, $(x-a_1)$, $\ldots$,
$(x-a_n)$, are irreducible in $\Int(\intz)$.
\endproof

\profess{Remark} Thanks to Alfred Geroldinger for pointing
this out: Theorem~\noblockmonoidthm\  implies that there does not 
exist a transfer-homomorphism from the multiplicative monoid
$(\Int(\Z)\setminus \{0\}, \cdot)$ to a block-monoid.
(For the definition of block-monoid and transfer-homomorphism
see {\rm \cite{GerHKNUF06}} Def.~2.5.5 and Def.~3.2.1, respectively.)

This is so because, in a block-monoid, the length of factorizations
of elements of the form $cd$ with $c$, $d$ irreducible, $c$ fixed,
is bounded by a constant depending only on $c$, cf.~{\rm \cite{GerHKNUF06}},
Lemma~6.4.4. 
More generally, applying {\rm \cite{GerHKNUF06}}, Lemma 3.2.2, one sees that
every monoid that admits a transfer-homomorphism to a block-monoid has
this property, in marked contrast to Theorem \noblockmonoidthm.
\endprofess
\bigskip

\line{\bf References\hfil}%
\nocite{*}
\bibliography{length}

\begin{thebibliography}{1}

\bibitem{CahChaEIVP95}
{\sc P.-J. Cahen and J.-L. Chabert}, {\em Elasticity for integral-valued
  polynomials}, J. Pure Appl. Algebra 103 (1995), 303--311.

\bibitem{CaCh97ivp}
{\sc P.-J. Cahen and J.-L. Chabert}, {\em Integer-valued polynomials}, vol.~48
  of Mathematical Surveys and Monographs, Amer.~Math.~Soc., 1997.

\bibitem{ChapMcCFE05}
{\sc S.~T. Chapman and B.~A. McClain}, {\em Irreducible polynomials and full
  elasticity in rings of integer-valued polynomials}, J. Algebra 293 (2005),
  595--610.

\bibitem{FreFriNUF11}
{\sc {\hbox{Ch}}.~Frei and S.~Frisch}, {\em Non-unique factorization of
  polynomials over residue class rings of the integers}, Comm. Algebra 39
  (2011), 1482--1490.

\bibitem{GerHKNUF06}
{\sc A.~Geroldinger and F.~Halter-Koch}, {\em Non-unique factorizations},
  vol.~278 of Pure and Appl. Math., Chapman \& Hall/CRC, Boca Raton, FL, 2006.

\bibitem{KaiFKIKG99}
{\sc F.~Kainrath}, {\em Factorization in {K}rull monoids with infinite class
  group}, Colloq. Math. 80 (1999), 23--30.

\end{thebibliography}
\bibliographystyle{siamese}
\bigskip
\goodbreak
\line{
\vbox{
\hbox{Institut f\"ur Mathematik A\hfil}
\hbox{Technische Universit\"at Graz\hfil}
\hbox{Steyrergasse 30\hfil}
\hbox{A-8010 Graz, Austria\hfil}
\tt
\hbox{frisch@tugraz.at\hfil}
}
\hfill
}

\bye